\documentclass[11pt]{amsart}
\usepackage{amsmath}
\usepackage{amsxtra}
\usepackage{amscd}
\usepackage{amsthm}
\usepackage{amsfonts}
\usepackage{amssymb}
\usepackage{eucal}

\newcommand{\nc}{\newcommand}

\numberwithin{equation}{section}

\swapnumbers

\newtheorem{thm}[subsection]{Theorem}
\newtheorem{prop}[subsection]{Proposition}
\newtheorem{lem}[subsection]{Lemma}

\newtheorem{rmk}[subsection]{Remark}

\nc{\ssec}{\subsection}
\nc{\sssec}{\subsubsection}

\nc{\on}{\operatorname}
\nc{\dmo}{\DeclareMathOperator}
\nc{\ZZ}{{\mathbb Z}}
\nc{\NN}{{\mathbb N}}
\nc{\FF}{{\mathbb F}}
\nc{\CC}{{\mathbb C}}

\nc{\F}{{\mathcal F}}
\nc{\W}{{\mathcal W}}
\nc{\J}{{\mathcal J}}

\nc{\Wb}{\overline{\W}}
\nc{\M}{{\mathcal M}}
\nc{\Y}{{\mathcal Y}}
\nc{\E}{{\mathcal E}}
\nc{\D}{{\mathcal D}}

\nc{\Z}{{\mathbb Z}}
\nc{\R}{{\mathbb R}}
\nc{\Q}{{\mathbb Q}}
\nc{\C}{{\mathbb C}}
\nc{\T}{{\mathbb T}}
\nc{\N}{{\mathbb N}}
\nc{\A}{{\mathbb A}}
\nc{\G}{{\mathbb G}}

\nc{\Ad}{\operatorname{Ad}}
\nc{\ad}{\operatorname{ad}}

\nc{\iso}{\buildrel \sim \over \longrightarrow}

\nc{\bi}{\bibitem}

\title{A Geometric Jacquet Functor}

\author{M. Emerton}

\address{Department of Mathematics, Northwestern University, Evanston, IL
60208, USA}

\author{D. Nadler}

\address{Department of Mathematics, University of Chicago, Chicago, IL
60637, USA}

\author{K. Vilonen}

\address{Department of Mathematics, Northwestern University, Evanston, IL
60208, USA}

\thanks{
The authors would like to acknowledge the partial support of
NSF grants DMS-0241562 and DMS-0401545 (M.E.),
an NSF postdoctoral fellowship (D.N.),
and NSF grant DMS-0105256 and DARPA grants MDA972-03-1-0030 and
HR0011-04-1-0031 (K.V.).
}

\begin{document}

\maketitle

%\address
%Northwestern University \newline\indent
%Department of Mathematics\newline\indent
%2033 Sheridan Rd.\newline\indent
%Evanston, IL 60208, USA
%\smallskip
%University of Chicago \newline\indent
%Department of Mathematics\newline\indent
%5734 S. University \newline\indent
%Chicago, IL 60637, USA
%\smallskip
%Northwestern University \newline\indent
%Department of Mathematics\newline\indent
%2033 Sheridan Rd.\newline\indent
%Evanston, IL 60208, USA
%\endaddress
%
%\email
%emerton\@math.northwestern.edu
%\smallskip
%nadler\@math.uchicago.edu
%\smallskip
%vilonen\@math.northwestern.edu
%\endemail

\section{Introduction}

In the paper \cite{Beilinson-Bernstein}, Beilinson and Bernstein
used the method of localisation to give a new proof and 
generalisation of Casselman's subrepresentation theorem. The key point is
to interpret $\mathfrak n$-homology in geometric terms.
The object of this note is to go one step further
and describe the Jacquet module functor on Harish-Chandra
modules via geometry.

Let $G_{\R}$ be a real reductive linear algebraic group,
and let $K_{\R}$ be a maximal compact subgroup of $G_{\R}$.
We use lower-case gothic letters to denote the corresponding Lie algebras,
and omit the subscript ``$\R$'' to denote complexifications.
Thus $(\mathfrak g, K)$ denotes the Harish-Chandra pair corresponding
to $G_{\R}$. 

Let $\mathfrak h$ be the universal Cartan
of $\mathfrak g$, that is
$\mathfrak h=\mathfrak b/ [\mathfrak b, \mathfrak b]$ where
$\mathfrak b$ is any Borel of
$\mathfrak g$. 
We equip $\mathfrak h$ with the usual choice of
positive roots by declaring
the roots of $\mathfrak b$ to be negative. 
We write $\rho\in{\mathfrak h}^*$ for half the sum of the positive roots.
To any
$\lambda\in {\mathfrak h}^*$ we associate a character~$\chi_\lambda$ of
the centre $Z(\mathfrak g)$ of the universal enveloping algebra
$U(\mathfrak g)$ via the Harish-Chandra homomorphism. 
Under this correspondence, the element $\rho\in\mathfrak h^*$
%half the sum of the positive roots,
corresponds to the trivial character $\chi_\rho$. For the
rest of this paper, we work with $\lambda\in\mathfrak h^*$ that is dominant,
i.e. $\check{\alpha}(\lambda) \not\in \{-1,-2,\ldots\}$ for any
positive coroot $\check{\alpha}$.

A Harish-Chandra module with infinitesimal character
$\lambda$ is, by definition, a $(\mathfrak g, K)$-module which is finitely
generated over $U(\mathfrak g)$ and on which %the centre 
$Z(\mathfrak g)$
acts via the character $\chi_\lambda$. We can also view Harish-Chandra
modules with infinitesimal character $\lambda$ as modules over the ring
$U_\lambda$ %=U(\mathfrak g)/I_\lambda$ 
which is the quotient of $U(\mathfrak g)$ by 
%where $I_\lambda$ is 
the two-sided ideal
generated by $\{z-\chi_\lambda(z) \mid z \in Z(\mathfrak g)\}$. In light
of this,
we will sometimes refer to such Harish-Chandra modules 
simply as $(U_\lambda, K)$-modules.

Let $X$ be the flag manifold of $\mathfrak g$, and let $\mathcal
D_\lambda$ be the sheaf of twisted differential operators with twist
$\lambda$. By a $(\mathcal D_\lambda, K)$-module we mean a coherent
$\mathcal D_\lambda$-module which is $K$-equivariant. Such a $\mathcal
D_\lambda$-module is, by necessity, regular holonomic since $K$ acts on $X$
with finitely many orbits. According to Beilinson-Bernstein
\cite{Beilinson-Bernstein2}, we have $\Gamma(X,\mathcal D_\lambda)=U_\lambda$, and
the global sections functor
\begin{equation}\label{BB}
\Gamma:\{(\mathcal D_\lambda, K)\text{-modules}\} \longrightarrow
\{(U_\lambda, K)\text{-modules}\}
\end{equation}
is exact and essentially surjective.
A section of the functor $\Gamma$
is given by the localisation functor,
which takes  
a $(U_\lambda, K)$-module $M$
to the
$(\mathcal D_\lambda, K)$-module
$\mathcal D_\lambda\otimes_{U_\lambda} M.$
The localisation functor is an equivalence if $\lambda$ is regular.

Let $P_{\R}$ be a minimal parabolic subgroup of $G_{\R}$,
let $\mathfrak n_\mathbb R$ be the nilpotent radical of $\mathfrak p_\mathbb R$, and
%We have $\mathfrak g_\mathbb R = \mathfrak p_\mathbb R \oplus
%\bar{\mathfrak n}_\mathbb R$ as $\mathbb R$-vector spaces, 
%where 
let $\bar{\mathfrak n}_\mathbb R$ be the
nilpotent subalgebra of $\mathfrak g_\R$
opposite to ${\mathfrak n}_\mathbb R$.
If $A_{\R}$ denotes the maximal split subtorus of the centre
of a Levi factor of $P_{\R}$, then $G_{\R}$ admits the 
Iwasawa decomposition $G_{\R} = K_{\R} A_{\R} N_{\R}.$

Given a
$(U_\lambda, K)$-module~$M$, we define 
the Jacquet module of $M$
by the formula
\begin{equation}
J(M)\ = \ \text{$\mathfrak n$-finite vectors in $\hat M$,}
\end{equation}
where $\hat M$ 
is the $\bar{\mathfrak n}$-adic completion $\varprojlim M/\bar{\mathfrak n}^k M$
of $M$.
% that is $\hat M = \varprojlim M/\bar{\mathfrak n}^kM$. 
We call the
functor $J$ which takes $M$ to $J(M)$ the Jacquet functor. As we will
explain in the next section, this functor, which is covariant and
preserves the infinitesimal character, is dual to the usual Jacquet
functor. The module $J(M)$ is a $(U_{\lambda},N)$-module (in the
sense that it is a Harish-Chandra module for the pair
$(\mathfrak g, N)$ with infinitesimal character $\lambda$).

To define the geometric Jacquet functor, we proceed as follows.
Let $\nu: \G_m \rightarrow A$ be a cocharacter of $G$ which is positive
on the roots in $\mathfrak
n$.  By composing  $\nu$ with the left action of $G$ on $X$, we obtain
an action of $\G_m$ on $X$ with action map
\begin{equation}\label{action}
a: \G_m \times X \rightarrow X.
\end{equation}
Consider the diagram
\begin{equation}
\G_m \times X \xrightarrow j \A^1 \times X \xleftarrow i \{ 0 \} \times X
\iso X
\end{equation}
where the maps $j$ and $i$ are the obvious inclusions.
To a $(\mathcal D_\lambda,K)$-module $\mathcal M$, we associate a
$\mathcal D_\lambda$-module $\Psi(\mathcal M)$ by taking the nearby
cycles of  $j_*a^*\mathcal M$  along $\{0\} \times X \iso X.$

Before stating our main theorem,
let us note that the localisation theory discussed above
for $(U_{\lambda}, K)$-modules
applies equally well to $(U_\lambda, N)$-modules.
%(The point is that $N$ also has only finitely many orbits on $X$.)

\begin{thm}\label{theorem}
(i) The localisation functor $M \mapsto
\mathcal D_\lambda\otimes_{U_\lambda} M$
takes the Jacquet functor $J$ to the geometric Jacquet functor
$\Psi$, i.e.,
$$
\mathcal D_\lambda \otimes_{U_\lambda} J(M) =
\Psi(\mathcal D_\lambda \otimes_{U_\lambda} M).
$$

(ii)
The global sections functor $\Gamma$ takes the geometric Jacquet functor
$\Psi$ to the Jacquet functor $J$, i.e., 
$$\Gamma(\Psi(\mathcal M))=
J(\Gamma(\mathcal M)).$$
\end{thm}

\begin{rmk} {\em  While the proof that we give of the main theorem is 
entirely algebraic,
its motivation is quite geometric.  Section \ref{Example} below provides
a discussion of its geometric interpretation: via the deRham functor $\operatorname{DR}$ the category of
$(\mathcal D_\lambda,K)$-modules is equivalent to the category of $\lambda$-twisted $K$-equivariant
perverse sheaves on $X$, and via this equivalence, $\Psi$ may be viewed as the operation of nearby cycles
on these perverse sheaves. We also discuss explicitly the case of 
$G_\mathbb R =
\on{SL}_2(\R)$.}
\end{rmk}

\begin{rmk} {\em The functor $\Psi$ takes $(\mathcal D_\lambda,K)$-modules to  $(\mathcal
D_\lambda,N)$-modules. This follows immediately from the above theorem and the fact that by the
translation principle we may assume that
$\lambda$ is  regular. 
%Alternatively, in section \ref{Example} it is shown that 
%$\operatorname{DR}(\Psi(\mathcal M))$ is constructible with respect to $P$-orbits.
}
\end{rmk}

\begin{rmk} {\em In this paper we could replace the ``Iwasawa nilpotent
subalgebra" $\mathfrak n$ by any larger nilpotent subalgebra.}
\end{rmk}

\begin{rmk} {\em The problem of giving a geometric construction of
the Jacquet functor was also considered by Casian and Collingwood.
In their paper \cite{Casian-Collingwood1},
they describe a nearby cycles construction on Harish-Chandra sheaves
(with dominant integral infinitesimal character);
they refer to the induced functor on $K$-groups (taking into account weight filtrations)
as a mixed Jacquet functor.  They conjecture that
this functor induces the usual Jacquet functor on Harish-Chandra
modules.

The nearby cycles construction
of \cite{Casian-Collingwood1} coincides with the construction appearing
in the definition of
the functor $\Psi$.  (Loosely, they pass to the limit with respect to
one simple coroot at a time,
while in the construction of~$\Psi$, we pass to the limit with respect
to all simple coroots at once.)  Thus, Theorem \ref{theorem}  establishes
the conjecture of Casian and Collingwood.}
\end{rmk}

%%%%%%%%%%%%%%%%%%%%%%%%%%%%%%%%%%%%%%%%%%%

\section {Jacquet functors}\label{Jacquet}

In this section we recall some standard facts about the Jacquet functor
as it is usually defined, and discuss the relation between this functor
and the functor $J$ defined in the Introduction.
As a general reference for basic facts, we use \cite{Wallach}, Chapter 4.

If $M$ is a $U_\lambda$-module, then the Jacquet
module of $M$ with respect to $\bar{\mathfrak n}$ is defined by
the formula
\begin{equation}
\check J(M) \ = \ \text{$\bar{\mathfrak n}$-finite vectors in $M^*$},
\end{equation}
where $M^*$ is the algebraic dual of $M$.  Of course, the infinitesimal
character of  $\check{J}(M)$ is $-\chi_\lambda = \chi_{\lambda^*}$,
where $\lambda^* = -\lambda + 2\rho$.

\begin{lem} The functor $\check{J}$ is an exact contravariant
functor from $(U_\lambda,K)$-modules (or more generally,
$U_{\lambda}$-modules that are finitely generated over $U(\bar{\mathfrak n})$)
to $(U_{\lambda^*},\bar N)$-modules.
\end{lem}
\begin{proof} The Lemma of Osborne
guarantees that any $(U_{\lambda},K)$-module is finitely generated over
$U(\bar{\mathfrak n})$.  If $M$ is such a finitely generated module,
then $\check{J}(M)$
is finitely generated over $U(\mathfrak g)$, and
for any positive integer $k$, the space of
$\bar{\mathfrak n}^k$-torsion vectors in $\check{J}(M)$
is finite-dimensional.  (It is dual to the finite-dimensional
space $M/\bar{\mathfrak n}^k M$.)  Thus 
$\check{J}(M)$ is a Harish-Chandra module for $(\mathfrak g, \bar N)$.

One may regard $\check{J}(M)$ as the topological
dual of the $\bar{\mathfrak n}$-adic completion $\hat{M}$ of $M$
(equipped with its $\bar{\mathfrak n}$-adic topology).  The exactness
of $\check{J}$ is thus implied by the exactness of forming
$\bar{\mathfrak n}$-adic completions.  This in turn follows from
the fact that since $\bar{\mathfrak n}$ is nilpotent the usual
Artin-Rees Lemma holds for the non-commutative ring $U(\bar{\mathfrak
n})$.
\end{proof}

Recall our definition of the Jacquet module $J(M)$ as $\mathfrak n$-finite
vectors in the $\bar{\mathfrak n}$-adic completion $\hat M$ of $M$.

\begin{lem} The functor $J$ is an exact covariant
functor from $(U_\lambda,K)$-modules to $(U_\lambda,N)$-modules.
\end{lem}

\begin{proof} First note that one may view $\hat M$ as the algebraic
dual $\check{J}(M)^*$ of $\check{J}(M)$. The decomposition $\mathfrak g =
{\mathfrak n} \oplus \bar{\mathfrak p}$ allows us to write $U(\mathfrak g)
= U({\mathfrak n}) \otimes_{\C} U(\bar{\mathfrak p})$ as $\C$-vector
spaces.  If $M$ is a $(U_\lambda,K)$-module then 
for a sufficiently large choice of $k$,
$\check{J}(M)$ is
generated over $U(\mathfrak g)$ by 
the space $\check{M} \, [\bar{\mathfrak n}^k]$
of $\bar{\mathfrak n}^k$-torsion elements in $\check{M}$. Since $U(\bar{ \mathfrak p})$ leaves
$\check{M}\, [\bar{\mathfrak n}^k]$ invariant, we see that $\check{J}(M)$
is finitely generated over $U({\mathfrak n})$. The lemma is now a
consequence of the preceding one (with $\mathfrak n$
replaced by $\bar{\mathfrak n}$).
\end{proof}

We will now give a more concrete description of
$J(M)$ which will be useful for our purposes. Recall from the last section
that we have a cocharacter 
$\nu: \G_m \to G$ which when
paired with the roots of $\mathfrak g$ is positive precisely on the roots in $\mathfrak n$. Let
us write  $h$ for the semisimple element in
$\mathfrak g$ given by $d\nu(t\partial_t)$. 
(%We think of the Lie algebra $\mathfrak g$ and that of $\mathbb G_m$ 
%as consisting of left-invariant vector fields. 
For any coordinate $t$ on $\mathbb G_m$,
the Euler vector field $t\partial_t$ gives %a bi-invariant vector field on $\mathbb G_m$, and so 
an element of the Lie algebra of $\mathbb G_m$
independent of the choice of coordinate.)
Note that the centraliser of $h$ is
precisely the Levi factor $\mathfrak l$ of $\mathfrak p$, and
%by construction, 
the weights of $h$ acting on $\mathfrak g$ are integral.

For each natural number $k$, the quotient $M/\bar{\mathfrak n}^kM$
is a finite-dimensional vector space on which $h$ acts,
and so may be written as a direct sum of generalised $h$-eigenspaces.
The surjection $M/\bar{\mathfrak n}^{k+1}M \rightarrow M/\bar{\mathfrak
n}^kM$ is $h$-equivariant, and so induces surjections of the
corresponding generalised $h$-eigenspaces.

As usual, we write $\hat{M}$ for the $\bar{\mathfrak n}$-adic completion
$\varprojlim M/\bar{\mathfrak n}^k M$ of $M$. 

\begin{lem} \label{eigenspaces}
If $M$ is a $(U_\lambda,K)$-module then the generalised
$h$-eigenspaces of %its $\bar{\mathfrak n}$-adic completion 
$\hat M$ are finite dimensional. For any integer $m$
the sum of generalised $h$-eigenspaces with eigenvalue greater than $m$
is finite dimensional. Furthermore,  $\hat M$ is the direct product of its
generalised $h$-eigenspaces.
\end{lem}

\begin{proof}
As the space $M/\bar{\mathfrak n}M$ is finite dimensional, it gives
rise to a finite set $S$ of generalised $h$-eigenvalues. The elements
in $\bar{\mathfrak n}^kM/\bar{\mathfrak n}^{k+1}M$ can be obtained by
multiplying $M/\bar{\mathfrak n}M$ by $\bar{\mathfrak n}^k$, and so each
eigenvalue of $h$ on $\bar{\mathfrak n}^kM/\bar{\mathfrak n}^{k+1}M$ is a
sum of an element of the set $S$ with an integer less than or equal to $-k$. This
gives our conclusion.
\end{proof}

\begin{prop}  \label{Jacquet and eigenspaces}
If $M$ is $(U_\lambda,K)$-module, then
the $(U_\lambda,N)$-module $J(M)$ is naturally isomorphic to
the direct sum of the generalised $h$-eigenspaces in $\hat{M}$.
\end{prop}

\begin{proof}
To prove the proposition, we must show that a vector in
$\hat{M}$ is ${\mathfrak n}$-finite if and only if it
is a sum of generalised $h$-eigenvectors.  Clearly any
${\mathfrak n}$-finite vector is such a sum since
it lies in a finite-dimensional $h$-invariant subspace
of $\hat{M}$.  Conversely, by the previous lemma,
the generalised $h$-eigenvalues of $\hat{M}$
are bounded above, and so
any sum of generalised $h$-eigenvectors is ${\mathfrak n}$-finite.
\end{proof}

For $\alpha \in \C$, let $\hat{M}_{\alpha}$ denote
the $\alpha$-generalised eigenspace of $h$ acting on $\hat{M}$.
We define an increasing $\C$-filtration on $\hat{M}$ by the formula
\begin{equation}\label{filtration}
F_{\alpha} (\hat{M}) = \prod_{\beta \leq \alpha} \hat{M}_{\beta}.
\end{equation}
Here and in what follows, 
for $\alpha,\beta\in\C$,
we write $\beta\leq\alpha$ to mean $\alpha-\beta$ is a non-negative
integer.
Pulling back this filtration via the injection $M \rightarrow \hat{M},$
we obtain a filtration on $M$.  The induced map on associated graded
modules
\begin{equation}
\on{Gr}^F_{\bullet} M \rightarrow \on{Gr}^F_{\bullet} \hat{M}
\end{equation}
is an isomorphism.  Indeed the proof of Lemma \ref{eigenspaces}
shows that the filtration just constructed
is cofinal with the $\bar{\mathfrak n}$-adic filtration
on $M$, and thus that $\hat{M}$ is the completion of $M$
with respect to this filtration.

Similarly, we
grade $\bar{\mathfrak n},$ and hence its enveloping
algebra $U(\bar{\mathfrak n})$, according to the eigenvalues
of the adjoint action of $h$.  Let $U(\bar{\mathfrak n})_\beta$
denote the graded piece on which $h$ acts with eigenvalue $\beta$.
Note that for $U(\bar{\mathfrak n})_{\beta}$ to be non-zero,
the eigenvalue $\beta$ must be a non-positive integer.
We define a filtration on $U(\bar{\mathfrak n})$ by the formula
\begin{equation}
F_{-k}(U(\bar{\mathfrak n})) = \bigoplus_{\beta \leq -k} 
U(\bar{\mathfrak n})_{\beta}.
\end{equation}
The $U(\bar{\mathfrak n})$-module
structure on $M$ is compatible with the filtrations on $U(\bar{\mathfrak n})$
and $M$.

\begin{lem} \label{Standard argument}  Let $\alpha$ be a complex number.

(i) For any $k\geq 0,$ the quotient
$F_{\alpha-k}(M)/(F_{-k}(U(\bar{\mathfrak n})) F_{\alpha}(M))$ is 
finite-dimensional.

(ii) If $\alpha \in \C$ is sufficiently small, then for any integer $k \geq 0$,
there is an equality
$F_{\alpha - k}(M) = F_{-k}(U(\bar{\mathfrak n})) F_{\alpha} (M)$
(and so the quotient considered in~(i) actually vanishes).
\end{lem}

\begin{proof}
The Lemma of Osborne guarantees that $M$ is finitely
generated over $U(\bar{\mathfrak n})$.  Also, the filtrations
constructed on each of $U(\bar{\mathfrak n})$ and $M$
are cofinal with the $\bar{\mathfrak n}$-adic filtrations.
Thus the lemma follows from an easy variant of the Artin-Rees Lemma (applied to
the filtrations $F_{\bullet}$ rather than the $\bar{\mathfrak n}$-adic
filtrations).
\end{proof}

%%%%%%%%%%%%%%%%%%%%%%%%%%%%%%%%%%%%%%%%%%%

\section{Nearby cycles, $V$-filtration, and the geometric Jacquet
functor}\label{geometric Jacquet}

In this section, we recall the formalism of the $V$-filtration
and the way it is used to define the  nearby cycles functor on $\mathcal
D$-modules. For a reference, see \cite{Kashiwara}.

As before, we write $X$ for the flag manifold of $\mathfrak g$, and 
$\mathcal D_\lambda$ for the sheaf of twisted differential operators with twist $\lambda$.
We denote by $t$ a coordinate
on $\A^1$, 
and consider the variety $\A^1 \times X$ and its subvariety $\{0\}\times X \simeq
X$. We write $\tilde{\mathcal D}_\lambda$ for the sheaf of
$\lambda$-twisted differential operators on $\A^1 \times X$. (Of course,
the twisting only occurs along the $X$-factor.) 
The $V$-filtration on $\tilde{\mathcal D}_\lambda$
is the filtration where we declare that the variable $t$ is of degree one,
$\partial_t$ is of degree $-1$, and everything in $\mathcal D_\lambda$ is
of degree zero.

For a regular holonomic $\tilde{\mathcal D}_\lambda$-module
$\mathcal M$ on $\A^1 \times X$, 
the $V$-filtration $V^\alpha \mathcal M$ 
is a decreasing filtration on
$\mathcal M$,
indexed by $\alpha\in\mathbb C$,
compatible with the
$V$-filtration on $\tilde{\mathcal D}_\lambda$.
%For $\alpha\in \mathbb C$,
%we write
%%denote by $\operatorname{Gr}_V^\alpha \mathcal M$ the associated graded
%\begin{equation}
%\operatorname{Gr}_V^\alpha \mathcal M \ = \ V^\alpha(\mathcal M)/
%V^{>\alpha}(\mathcal M).
%\end{equation}
The  V-filtration on $\mathcal M$ is uniquely determined
by the following properties:
\begin{equation}\label{V1}
\begin{split}
\text{there exist elements } & \text{$u_1, \dots, u_k\in \mathcal
M$ and
$\mu_1,\dots,\mu_k\in \C$ such that}
\\
  &V^\alpha(\mathcal M) = \sum _{\alpha\leq
m_i+\mu_i} V^{m_i}(\tilde{\mathcal D}_\lambda)u_i,
\end{split}
\end{equation}
and
\begin{equation}\label{V2}
\text{%for $\alpha\in \mathbb C$,
$\operatorname{Gr}_V^\alpha \mathcal M$ is a generalised
eigenspace of $t\partial_t$ of eigenvalue $\alpha$.}
\end{equation}
Here for $\alpha\in \mathbb C$,
we have written $\operatorname{Gr}_V^\alpha$ for the quotient
$
V^\alpha(\mathcal M)/V^{>\alpha}(\mathcal M).
$

The nearby cycles of $\mathcal M$ along $\{0\}\times X\simeq X$ are given by
\begin{equation}
\psi(\mathcal M) = \operatorname{Gr}_V^0 \mathcal
M\oplus\bigoplus_{\substack{\alpha\in
\C/\Z\\ \alpha\notin\Z}}\operatorname{Gr}_V^\alpha
\mathcal M.
\end{equation}

We are now ready to define the geometric Jacquet functor. 
Recall from the first section that via the cocharacter $\nu:\mathbb G_m\to G$, 
we have defined %in \eqref{action} 
a $\G_m$-action on $X$. This gives us the diagram
\begin{equation}
X \xleftarrow{a}\G_m \times X \xrightarrow j \A^1 \times X \xleftarrow i
\{ 0\}\times X\iso X
\end{equation}
where $a$ is the action map, and the maps $j$ and $i$ are the obvious inclusions.
We define the geometric Jacquet functor $\Psi$ by the formula
\begin{equation}
\Psi(\mathcal M) \ = \ \psi (j_*a^*\mathcal M),\qquad\text{for a
$(\mathcal D_\lambda,K)$-module $\mathcal M$}.
\end{equation}

%%%%%%%%%%%%%%%%%%%%%%%%%%%%%%%%%%%%%%%%%%%%%%%

\section{Proof of the main theorem}\label{Proof of the main theorem}

In this section, we give a proof of Theorem \ref{theorem}. To
begin with, suppose that $\mathcal M$ is a $(\mathcal D_\lambda,K)$-module
arising as the localisation of a $(U_\lambda,K)$-module $M$,
so that $\mathcal M = \D_{\lambda}\otimes_{U_{\lambda}} M.$
%Since $M$ is the global sections of $\mathcal M$,
Our aim is to explicitly write down
the $V$-filtration on $j_*a^*\mathcal M$ %may be calculated 
on the level of its global sections $\tilde M$.
%Let $\tilde M$ denote the global sections 
%$\Gamma(\A^1\times X,j_*a^*\mathcal M)$. 
We will make use of the diagram
\begin{equation}\label{pullback}
\G_m \times X \xrightarrow {\tilde a} \G_m \times X \xrightarrow p X,
\end{equation}
where $p$ is the projection, and the map $\tilde a$ is given by $\tilde
a(g,x)= (g,gx)$ so that $a=\tilde a \circ p$. 

To describe $\tilde M$ as a vector space, 
we use the coordinates
in the middle copy of $\G_m \times X$ in diagram \eqref{pullback}, that is
we take the global sections of $p^*\mathcal M$. This gives us
\begin{equation}
\tilde M \ = \ \Gamma(\A^1\times X,j_*a^*\mathcal M)
\ = \ \C[t,t^{-1}]\otimes_\C M.
\end{equation}

To describe the action of 
\begin{equation}
\tilde U_\lambda = \Gamma(\G_m \times
X,\tilde{\mathcal D}_\lambda)=\C[t,t^{-1},t\partial_t]\otimes_\C
U_\lambda,
\end{equation}
we use the coordinates in
the first copy of $\G_m \times X$ in diagram \eqref{pullback}.
It suffices to
write down the action %on $\tilde M$ %, which we denote by `$\cdot$',
of the global vector field $t\partial_t$ and the global vector fields given by elements of
$\mathfrak g$. 
To this end, let us recall some of our previous notation: we have the
cocharacter $\nu : \G_m \to G$, and we write $h$ for the semisimple
element in $\mathfrak g$ given by $d\nu(t\partial_t)$. 
First, under the map $\tilde a$ the global vector
field $(t\partial_t,0)$ is mapped to the global vector field
$(t\partial_t, h)$. Hence,
\begin{equation}\label{twist1}
t\partial_t\cdot (f(t)m)\ = \ (t\partial_t f(t))m +
f(t)(hm).
\end{equation}
Second, under the map $\tilde a$ the global vector field $(0,v)$
given
by an element $v\in\mathfrak g$ is mapped to the global vector field
$(0,\Ad_{\nu^{-1}(t)} v)$. Hence,
\begin{equation}\label{twist2}
v\cdot (f(t)m)\ =  \ f(t) ((\Ad_{\nu^{-1}(t)} v) m).
\end{equation}
In particular, if $v$ is in the $\beta$ root space $\mathfrak g_\beta$, then
\begin{equation}\label{twist3}
v\cdot m\ = \ t^{-<\nu,\beta>} vm.
\end{equation}
(The inverse of the adjoint action arises for the following reason.
For a vector $v\in\mathfrak g$ thought of as a vector field on $X$,
the bracket $[h,v]$ may be calculated by pulling back~$v$ along the
integral curves given by the action of $\nu(t)$. Thus
the derivative of the pushforward of $v$ under the action of $\nu(t)$ is
the negative of $[h,v]$.)
% sign in \ref{twist2} and \ref{twist3}...)

To write down the $V$-filtration on $\tilde M$, 
let us recall the discussion of Section \ref{Jacquet}. There we
constructed an increasing $\C$-filtration $F_{\bullet} (M)$ on $M$.
We now argue that the $V$-filtration on $\tilde M$
may be written in terms of this filtration as follows
\begin{equation}
V^\alpha (\tilde M) \ = \ \bigoplus_{k\in \Z} t^k F_{-\alpha+k}(M) \subset
\bigoplus_{k\in \Z} t^k M=\tilde M.
\end{equation}
To show that this is the $V$-filtration we must verify conditions
\eqref{V1} and \eqref{V2}. The fact that \eqref{V2} holds is clear from
the construction (taking into account equation \ref{twist1}).
To show that \eqref{V1} holds, we show that
\begin{equation}\label{V11}
V^{\alpha+m}(\tilde M) \ = \ V^m(\tilde U_\lambda)V^{\alpha}(\tilde M),
\qquad \text{for all $\alpha$ and $m$},
\end{equation}
and that
\begin{equation}\label{V12}
\text{each $V^{\alpha}(\tilde M)$ is finitely generated over $V^0(\tilde
U_\lambda)$}.
\end{equation}
The first statement \eqref{V11} is clear since multiplication by $t$ clearly
induces an isomorphism between $V^{\alpha}(\tilde M)$ and
$V^{\alpha+1}(\tilde M)$. Therefore it remains to argue \eqref{V12}. We will
in fact show that $V^{\alpha}(\tilde M)$ is finitely generated over
$\C[t]\otimes_{\C}U(\bar{\mathfrak n})$. By Lemmas
\ref{eigenspaces} and \ref{Standard argument}, we have
\begin{equation}
\text{ $F_{-\alpha + k}(M) =M$, for sufficiently large $k$},
\end{equation}
\begin{equation}\label{finite dimensionality}
\text{ $F_{-\alpha -k}(M)/(F_{-k}(U(\bar{\mathfrak n})) F_{\alpha}(M))$
is finite-dimensional, for all $k \geq 0$},
\end{equation}
and
\begin{equation}\label{stability}
\text{ $F_{\alpha - k - l}(M) = F_{-l}(U(\bar{\mathfrak n})) 
F_{\alpha - k} (M)$,
for sufficiently large $k$ and all $l\geq 0$}.
\end{equation}
These three statements together
(taking into account equation \ref{twist3})
imply that $V^{\alpha}(\tilde M)$ is
finitely generated over $\C[t]\otimes_{\C}U(\bar{\mathfrak n})$, and 
hence that it is
finitely generated over $V^0(\tilde U_\lambda)=\C[t,t\partial_t]\otimes_{\C}
U_\lambda$.

Now by Section \ref{geometric Jacquet}, we have
\begin{equation}
\Psi(M) \ = \ \text{direct sum of generalised $h$-eigenspaces
in $\hat M$}.
\end{equation}
Part~(i) of Theorem \ref{theorem} thus follows from Proposition \ref{Jacquet and
eigenspaces}.

We turn to proving part~(ii) of the theorem.  Thus we suppose
that $\mathcal M$ is an arbitrary
$(\mathcal D_\lambda,K)$-module.
If we write $M = \Gamma(\mathcal M),$ then there is a
canonical map $\mathcal D_\lambda \otimes_{U_\lambda} M \rightarrow
\mathcal M,$ whose kernel and cokernel are both taken to zero
by the exact functor
$\Gamma$.  Taking into account
the fact that the passage to nearby cycles is exact, as well as part~(i)
of the theorem,
we see that part~(ii) follows from the next statement:
\begin{equation}
\label{vanishing}
\text{if } \Gamma(\mathcal M) = 0, \text{ then }
\Gamma (\Psi(\mathcal M)) = 0.
\end{equation}
To prove \eqref{vanishing}, we will 
argue as follows. First, it 
suffices to prove
it in the case when $\mathcal M$ is irreducible.
Then, we will appeal to the results of \cite{Kashiwara-singular}, \S 8.
The final remark of \cite{Kashiwara-singular}, p.~103,
shows that we may find a Weyl group
element $w$ so that $w(\lambda)$ is anti-dominant, and satisfies
condition~(8.3.2).  Applying the intertwining functor $I_w$ 
induces an equivalence of categories between the category
of $(\D_{\lambda},K)$-modules and the category of
$(\D_{w(\lambda)},K)$-modules, which is easily verified
to be compatible with the functor $\Psi$ on each of these categories.
Thus in verifying~\ref{vanishing} we may replace $\lambda$ by
$w(\lambda)$.  
The remark preceding the statement of \cite{Kashiwara-singular},
Prop.~8.2.1 shows that we may also replace $G$ by its simply-connected
cover, and thus assume that the sheaf $\mathcal O(\rho)$ is a $G$-equivariant
line-bundle on $X$.
Proposition~8.2.1 and Theorem~8.3.1 of \cite{Kashiwara-singular} now show
there exists a simple root $\alpha$ such that $\check\alpha(\lambda)=0$ and 
$\mathcal M = \mathcal O(-\rho)\otimes p^* \mathcal M'$,  where
$p:X\to X_\alpha$ is the projection from the full flag manifold $X$ 
to the partial flag manifold
$X_\alpha$ of parabolics of the type associated to the simple root 
$\alpha$, and $\mathcal M'$ is a $\D_{\lambda + \rho}$-module on $X_{\alpha}$.
The functor $\Psi$ has
its analogue $\Psi_\alpha$ defined on $\mathcal D_{\lambda+\rho}$-modules on 
$X_\alpha$, compatible with
$\Psi$. Thus we conclude
\begin{multline}
\Gamma (\Psi(\mathcal M)) =
\Gamma (\Psi(\mathcal O(-\rho) \otimes p^* \mathcal M')) =
\Gamma(\mathcal O(-\rho) \otimes \Psi(p^* \mathcal M')) \\
= \Gamma (\mathcal O(-\rho) \otimes p^*\Psi_\alpha(\mathcal M')) =
\Gamma (p_*(\mathcal O(-\rho) \otimes p^*\Psi_\alpha(\mathcal M'))) = 0;
\end{multline}
here the second equality follows from the evident compatibility of the
functor $\Psi$ with twisting by $G$-equivariant line-bundles,
the third equality follows from the compatibility of $p^*$ with
$\Psi$ and $\Psi_{\alpha}$, and the
the last equality holds because, by \cite{Kashiwara-singular}, Prop.~8.2.1,
we see that $p_*(\mathcal O(-\rho)\otimes p^*(\mathcal M'))$
is zero for any $\mathcal 
D_\lambda$-module $\mathcal M'$ on $X_\alpha$, when $\check\alpha(\lambda)=0$.

%%%%%%%%%%%%%%%%%%%%%%%%%%%%%%%%%%%%%%%%%%%%%
%%%%%%%%%%%%%%%%%%%%%%%%%%%%%%%%%%%%%%%%%%%%%
%%%%%%%%%%%%%%%%%%%%%%%%%%%%%%%%%%%%%%%%%%%%%

\section{A geometric interpretation}\label{Example}

As mentioned in the introduction, the arguments which go into the proof
of our main theorem are primarily algebraic. 
In this section,
we informally describe the geometry those arguments formalise and also discuss explicitly the example
$G_\R=\on{SL}_2(\R)$.
The discussion is in the language of $\lambda$-twisted perverse sheaves rather than holonomic $\mathcal D_\lambda$-modules.

To review our notation, $\mathfrak p_\R$ is a minimal parabolic subalgebra 
of $\mathfrak g_\R$
with complexification $\mathfrak p$;
$\mathfrak n$ is the nilpotent radical
of $\mathfrak p$;
$P$ is the parabolic subgroup with Lie algebra $\mathfrak p$;
$\bar P$ is the parabolic subgroup opposite to $P$;
$L$ is the Levi subgroup $P\cap \bar P$; and $X$ is the flag 
manifold of $\mathfrak g$.

The basic ingredient in defining the geometric Jacquet functor is the action
\begin{equation}
a: \G_m \times X \rightarrow X
\end{equation}
defined as the composition of the action of $G$ on $X$
with a cocharacter $\nu: \G_m \rightarrow G$.
%We shall think of this action as providing a flow on $X$.
%%given by 
%%$$
%%x\mapsto a(z)\cdot x, \mbox{ for } z\in\G_m.
%%$$
Recall that we choose $\nu$ such that when
paired with the roots of $G$ it is positive precisely on the roots in
$\mathfrak n$.

The action $a$ defines two Morse stratifications of
$X$, one by ascending manifolds, one by descending manifolds,
which have the same critical manifolds.
The ascending manifolds are the $P$-orbits,
%$X_P^\alpha$,
the descending manifolds are the $\bar P$-orbits,
%$X_{\bar P}^\alpha$ ,
and the critical manifolds
%$X^\alpha$
are the $L$-orbits. %intersections of corresponding orbits.
% $X_{P}^\alpha\cap X_{\bar P}^\alpha$.
Indeed, each connected component $X^\alpha$
of the fixed points of $a$ is a single $L$-orbit, the $P$-orbit $X^\alpha_P$
through $X^\alpha$ satisfies
\begin{equation}
X^\alpha_P=\{x\in X|\lim_{z\to0} a(z)\cdot x\in X^\alpha\},
\end{equation}
and the $\bar P$-orbit $X^\alpha_{\bar P}$
through $X^\alpha$ satisfies
\begin{equation}
X^\alpha_{\bar P}=\{x\in X|\lim_{z\to \infty} a(z)\cdot x\in X^\alpha\}.
\end{equation}
%In other words, 

Given a $\lambda$-twisted $K$-equivariant perverse sheaf $\mathcal M$ on $X$
%$(\mathcal D_\lambda,K)$-module $\mathcal M$ on $X$, 
we think of its geometric Jacquet module $\Psi(\mathcal M)$ as the result of "flowing"
$\mathcal M$ along the trajectories of the action $a$ and "passing to 
the limit".
The limit $\Psi(\mathcal M)$
is constructible with respect to the stratification by the ascending 
manifolds $\{X^\alpha_P\}$,
and an ascending manifold $X^\alpha_P$ is in the support of $\Psi(\mathcal M)$
exactly when the corresponding
descending manifold $X^\alpha_{\bar P}$ intersects the support of $\mathcal M$.
The restriction of $\Psi(\mathcal M)$
to a critical manifold $X^\alpha$ may be calculated as the integral 
of $\mathcal M$
over the fibers of the corresponding descending manifold $X^\alpha_{\bar P}$.

To see this in a concrete example, consider
the real group $\on{SL}_2(\R)$, with maximal compact subgroup
$\on{SO}_2(\R)$, and the corresponding
Harish-Chandra pair $(\mathfrak{sl}_2(\C),\on{SO}_2(\C))$.
We may identify the flag manifold of $\mathfrak{sl}_2(\C)$ with
the complex projective line $\C\mathbb P^1=\C\cup \{\infty\} $
so that (i) the $\on{SO}_2(\C)$-orbits
are the points $i$ and $-i$,
and their complement $\C\mathbb P^1\setminus \{i, -i\}$,
and (ii) the point $0$ represents a Borel subalgebra
$\mathfrak b\subset\mathfrak{sl}_2(\C)$
which is the complexification of a minimal parabolic subalgebra
of $\mathfrak p_\R\subset \mathfrak{sl}_2(\R)$.

Fix a character $\lambda:\mathfrak b/[\mathfrak b,\mathfrak b]\to\C$.
Let $\delta_0$ be the irreducible $\lambda$-twisted perverse sheaf 
%delta function $\mathcal D_\lambda$-module 
on $\C\mathbb P^1$
supported at the point $0$,
and let $\mathcal N$ be the middle-extension to $\C\mathbb P^1$
of the irreducible $\lambda$-twisted local system
%$\mathcal D_\lambda$-module 
on $\C\mathbb P^1\setminus\{0\}$.
For $\lambda$ non-integral, there are
two irreducible $\lambda$-twisted $\on{SO}_2(\C)$-equivariant perverse sheaves 
%$(\mathcal D_\lambda, \on{SO}_2(\C))$-modules 
on $\C\mathbb P^1$:
the middle-extensions $\mathcal M_{1}$ and $\mathcal M_{-1}$
of the irreducible trivial and twisted $\on{SO}_2(\C)$-equivariant local systems on
the complement $\C\mathbb P^1\setminus\{i,-i\}$.
For $\lambda$ integral, there are in addition two other irreducible
$\lambda$-twisted $\on{SO}_2(\C)$-equivariant perverse sheaves 
%$(\mathcal D_\lambda, \on{SO}_2(\C))$-modules 
on $\C\mathbb P^1$:
the irreducible %$\lambda$-twisted 
$\on{SO}_2(\C)$-equivariant 
perverse sheaves $\delta_i$ and $\delta_{-i}$ supported at the 
points $i$ and~$-i$.

Choose the cocharacter $\nu:\C^\times\to \on{SL}_2(\C)$ so that
the induced action $a:\C^\times\times\C\mathbb P^1\to\C\mathbb P^1$
is the standard multiplication $a(z)\cdot p=z^{-2} p$.
Note that in the limit $z\to \infty$, the points $i$ and $-i$ flow to 
$0$, and this
completely describes what happens to the $\on{SO}_2(\C)$-orbits.
%Given a $(\mathcal D_\lambda, \on{SO}_2(\C))$-module
%on $\C\mathbb P^1$, we calculate its geometric Jacquet module
%by flowing it along the trajectories of the action
%$a$ and passing to the
%limit as $z\to \infty$.
With this picture in mind, we calculate the geometric Jacquet module
of the irreducible $\lambda$-twisted $\on{SO}_2(\C)$-equivariant perverse sheaves
%$(\mathcal D_\lambda, \on{SO}_2(\C))$-modules
on $\C\mathbb P^1$.
For $\lambda$ integral, we have
$$
\Psi(\delta_i)=\Psi(\delta_{-i})=\delta_0
\quad
\Psi(\mathcal M_1)=\mathcal N
\quad
\Psi(\mathcal M_{-1})=\mathcal L
$$
where $\mathcal L$ is the unique non-semisimple self-dual
perverse sheaf with three-step filtration $\mathcal L_0\subset\mathcal 
L_1\subset\mathcal L$
satisfying
$$
\mathcal L_0=\delta_0
\quad
\mathcal L_1/\mathcal L_0 =\mathcal N
\quad
\mathcal L/\mathcal L_1 =\delta_0.
$$
For $\lambda$ non-integral, we have
$$
\Psi(\mathcal M_1)=\Psi(\mathcal M_{-1})=\mathcal 
N\oplus\delta_0\oplus \delta_0.
$$
To help explain why the result for $\mathcal M_1$ and $\lambda$ integral
has different simple constituents
from the other perverse sheaves supported on the open orbit, note that only 
for $\mathcal M_1$
and $\lambda$ integral is the characteristic variety of such a perverse sheaf
the zero section alone.

%%%%%%%%%%%%%%%%%%%%%%%%%%%%%%%%%%%%%%%%%%%%%
%%%%%%%%%%%%%%%%%%%%%%%%%%%%%%%%%%%%%%%%%%%%%

\end{document}